\documentclass[12pt]{article}
\usepackage[dvips]{graphics}
\usepackage{amsthm,amsfonts}
\usepackage{fullpage}

\newtheorem{theorem}{Theorem}

\begin{document}
\title{A note on non-repetitive colourings of planar graphs}
\author{Narad Rampersad \\
School of Computer Science \\
University of Waterloo \\
Waterloo, ON, N2L 3G1 \\
CANADA \\
{\tt nrampersad@math.uwaterloo.ca}}
\date{\today}
\maketitle

\begin{abstract}
Alon \emph{et al.} introduced the concept of non-repetitive
colourings of graphs.  Here we address some questions regarding
non-repetitive colourings of planar graphs.  Specifically, we
show that the faces of any outerplanar map can be non-repetitively
coloured using at most five colours.  We also give some lower bounds
for the number of colours required to non-repetitively colour the
vertices of both outerplanar and planar graphs.
\end{abstract}

\section{Introduction}
A sequence $a$ is called \emph{non-repetitive} if $a$ contains no
identical adjacent blocks.  A vertex (resp. edge) colouring
of a graph $G$ is called \emph{non-repetitive} if for any open
path $P$ in $G$, the sequence of vertex (resp. edge) colours
along $P$ is non-repetitive.  If $G$ is a planar map, then a colouring
of the faces of $G$ is called non-repetitive if for any sequence
of distinct faces such that each consecutive pair of faces shares
an edge, the sequence of corresponding colours is non-repetitive.

Grytczuk \cite{Gry02} asked the following question:  is there a natural number
$k$ such that the faces of \emph{any} planar map can be non-repetitively
coloured using at most $k$ colours?  We answer this question in the
affirmative for all \emph{outerplanar} maps.

Alon \emph{et al.} \cite{AGHR02} asked a similar question:  is there a
natural number $k$ such that the vertices of \emph{any} planar graph can be
non-repetitively coloured using at most $k$ colours?  Here we give lower
bounds for $k$ for both \emph{outerplanar} and \emph{planar} graphs.

\section{Main results}
\begin{theorem}
\label{map}
If $G$ is an outerplanar map, then the faces of $G$ can be non-repetitively
coloured using at most five colours.
\end{theorem}

\begin{proof}
We assume that we have an outerplanar embedding of $G$.
It suffices to show that the vertices of the dual graph $G^*$ can
be coloured non-repetitively using at most five colours.  Consider
the weak dual $G^w$ formed by deleting the vertex of $G^*$ corresponding
to the outer face of $G$.  It is well known \cite{FGH74} that if $G$ is
outerplanar, then $G^w$ is a forest of trees.  Alon \emph{et al.}
\cite{AGHR02} mention that the vertices of any tree can be
non-repetitively coloured using at most four colours.  Hence, the weak
dual $G^w$ can be non-repetitively coloured with at most four colours.
Finally, add back the vertex initially deleted from $G^*$ and
colour it using a fifth colour.  The resulting 5-colouring of the dual
graph $G^*$ induces a non-repetitive colouring of the faces of the
outerplanar map $G$.
\end{proof}

\begin{theorem}
\label{outer}
There exists an outerplanar graph $G$ such that the vertices
of $G$ cannot be non-repetitively coloured using fewer than five
colours.
\end {theorem}

\begin{center}
\begin{figure}[h]
\includegraphics{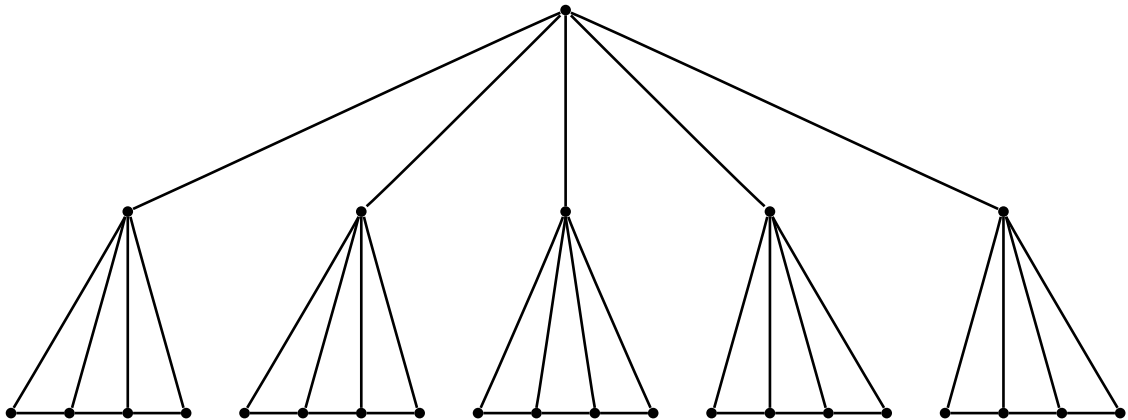}
\caption{\label{out_pic} Graph $G$ from Theorem~\ref{outer}}
\end{figure}
\end{center}

\begin{proof}
We will construct such a graph $G$.  We begin with the graph $P_4$,
\emph{i.e.} the graph consisting of a simple path on four vertices.
Since there are no non-repetitive binary sequences of length four,
we require at least three colours to non-repetitively colour $P_4$.
Next we add a vertex $v$, connecting it with an edge to each of the
vertices of $P_4$, thus forming the so-called \emph{fan graph} $F_4$.
Let us call the vertex $v$ the \emph{rivet} of the fan.
Since $v$ is connected to vertices of three different colours,
it is evident that $v$ must be coloured a fourth colour.  The graph
$G$ then consists of \emph{five} disjoint copies of $F_4$, with an
additional vertex $r$ connected to the rivet of each fan (see
Figure~\ref{out_pic}).  Clearly $G$ is an outerplanar graph.
If we assume that we only have four colours with which to work,
then by the pigeonhole principle, two rivets, say $v$ and $v'$,
must be coloured the same colour, say $x$.  The vertex $r$ cannot be
coloured $x$, so $r$ must be coloured with one of the three remaining
colours, say $y$.  However, the subgraph $P_4$ connected to the
rivet $v$ contains vertices coloured with three distinct colours different
from $x$.  Hence, we can always find a vertex $w$ such that the path
$wvrv'$ has colouring $yxyx$.  This is clearly a repetition, and so
we see that we need at least five colours to non-repetitively colour
$G$.
\end{proof}

\begin{theorem}
\label{planar}
There exists a planar graph $G$ such that the vertices of $G$
cannot be non-repetitively coloured using fewer than seven colours.
\end{theorem}

\begin{center}
\begin{figure}[h]
\includegraphics{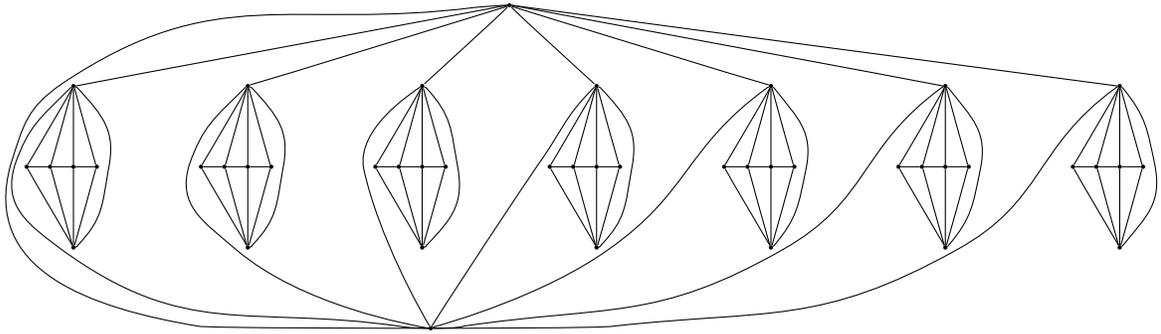}
\caption{\label{plan_pic} Graph $G$ from Theorem~\ref{planar}}
\end{figure}
\end{center}

\begin{proof}
The construction of $G$ is readily apparent from Figure~\ref{plan_pic}.
Let us label the two vertices of $G$ with degree eight $r$ and $s$.  Let
us call each of the connected components of the subgraph formed by
deleting $r$ and $s$ from $G$ a \emph{diamond}.  By reasoning similar
to that used in the proof of Theorem~\ref{outer}, we may conclude that
each diamond requires at least five colours for a non-repetitive
colouring.  Assume that we have a non-repetitive 6-colouring of $G$.
Now consider the seven vertices of $G$ with degree seven.  By the
pigeonhole principle, two of these vertices must have the same colour.
Let us call these two vertices $v$ and $v'$, and let us assume that they
are each coloured $x$.  Let us call each of the two diamonds containing
$v$ and $v'$ $D$ and $D'$ respectively. Suppose that $D$ and $D'$ are
each coloured using exactly five colours, but the five colours used
are not the same for each diamond.  In this case, between $D$ and $D'$
all six colours are used, and so for all choices $y$, $y \neq x$,
for the colour of $r$ we can always find a vertex $w$ in one of
$D$ or $D'$ such that the path $wvrv'$ has colouring $yxyx$.
Hence, it must be the case that $D$ and $D'$ are each
coloured using exactly the same colours.  If $D$ and $D'$ are each coloured
using all six colours, then again we can always find a vertex $w$ in one of
$D$ or $D'$ such that the path $wvrv'$ has colouring $yxyx$.  It is
therefore the case that $D$ and $D'$ are each coloured using exactly the
same five colours.  Thus we may colour $r$ using the colour that does not
appear in $D$ or $D'$; any other choice $y$, $y \neq x$, will allow us to
find a vertex $w$ in one of $D$ or $D'$ such that the path $wvrv'$ has
colouring $yxyx$.  However, by this same argument we find that $s$
must be coloured the same colour as $r$.  Since $r$ and $s$ are adjacent,
this is a contradiction, and so we have that $G$ cannot be non-repetitively
coloured using fewer than seven colours.
\end{proof}

\begin{thebibliography}{99}
\bibitem{AGHR02}
N. Alon, J. Grytczuk, M. Haluszczak, O. Riordan, ``Nonrepetitive
colorings of graphs'', \emph{Random Structures Algorithms} \textbf{21}
(2002), 336--346.
\bibitem{Gry02}
J. Grytczuk, ``Thue-like sequences and rainbow arithmetic progressions'',
\emph{Electron. J. Combin.} \textbf{9} (2002), \#R44.
\bibitem{FGH74}
H.J. Fleischner, D.P. Geller, F. Harary, ``Outerplanar graphs and weak
duals'', \emph{J. Indian Math. Soc.} \textbf{38} (1974), 215--219.
\end{thebibliography}
\end{document}